\documentclass[11pt]{article}
\newfont{\frak}{eufm10 scaled\magstep1}
\newfont{\sfrak}{eufm8 scaled\magstep1}
\newfont{\bbb}{msbm10 scaled\magstephalf}
\newfont{\sbbb}{msbm7 scaled\magstephalf}

\def\D{\Delta}
\def\C{\mbox{\bbb{C}}}
\def\R{\mbox{\bbb{R}}}
\def\Z{\mbox{\bbb{Z}}}
\def\cd{\C^d}
\def\rd{\R^d}

\def\vz{\underline{z}}
\def\ed{e_1,\ldots,e_d}
\def\ld{\lambda_1,\ldots,\lambda_d}
\def\xd{X_1,\ldots,X_d}
\def\d{\mbox{\frak d}}
\def\n{\mbox{\frak n}}
\def\ddu{\d^*}

\def\zd{(z_1,\cdots,z_d)}
\def\lorw{\longrightarrow}
\def\SC{\mbox{\sbbb{C}}}
\def\Dc{D_{\SC}}
\def\cdf{\C^d_F}
\def\zset{\Psi^{-1}(0)}

\def\cono{\stackrel{\circ}{C}}

\newtheorem{thm}{Theorem}[section]
\newtheorem{prop}[thm]{Proposition}
\newtheorem{lemma}[thm]{Lemma}
\newtheorem{defn}[thm]{Definition}
\newtheorem{remark}[thm]{Remark}
\newcommand{\proof}{\mbox{\textbf{ Proof.\ \ }}}

\title{\sc Nonrational, nonsimple convex polytopes in symplectic geometry}
\author{\sc Fiammetta Battaglia\thanks{Partially supported by
MIUR project
\textit{Propriet\`a Geometriche delle Variet\`a Reali e Complesse},
by  GNSAGA (CNR), and by EDGE (EC FP5 Contract no. HPRN-CT-2000-00101).}
and Elisa Prato}
\date{}
\begin{document}
\maketitle
\begin{abstract}
In this research announcement we associate to each convex
polytope, possibly nonrational and nonsimple, a family of compact
spaces that are stratified by quasifolds, i.e. the strata are
locally modelled by $\R^k$ modulo the action of a discrete,
possibly infinite, group. Each stratified space is endowed with a
symplectic structure and a moment mapping having the property that
its image gives the original polytope back. These spaces may be
viewed as a natural generalization of symplectic toric varieties
to the nonrational setting. We provide here the explicit
construction of these spaces, and a thorough description of the
stratification.
\end{abstract}

{\small 2000 \textit{Mathematics Subject Classification.} Primary: 53D05.
Secondary: 53D20, 32S60, 52B20.}


\parindent0pt

\section*{Introduction}

To each rational convex polytope it is possible to associate, by a
standard construction, a geometric object that is known as the
toric variety corresponding to the polytope. Is it possible to
associate a similar geometric object to a convex polytope that is
not rational?

In a number of recent papers this object has been referred to as
``nonexisting'' (e.g. \cite{bl} and, along the same lines,
\cite{bbfk,bbfk2}), but in fact it is shown by the authors in
\cite{p,cx} that for convex polytopes that are simple such an
object (and in fact a whole family of such objects) exists; it is
an example of a space, known as quasifold, that is locally modelled
by $\R^k$ modulo the action of a discrete, possibly infinite, group
and represents a natural generalization of a toric variety.

In the present announcement we consider the problem from the
symplectic viewpoint, in the case of convex polytopes that are no
longer simple.

We show that given an $n$-dimensional vector space, $\d$, and a
convex polytope $\D\subset\d^*$, there is a family of compact
spaces that are stratified by symplectic quasifolds. Each space
$M$ of the family admits the continuous action of an
$n$-dimensional quasitorus $D$ and a continuous mapping
$\Phi\,\colon\,M\lorw\d^*$ such that $\Phi(M)=\D$. The restriction
of the $D$-action to each stratum is smooth and Hamiltonian, with
moment mapping given by the restriction of $\Phi$. These
stratified spaces are a natural generalization to the quasifold
setting of the notion of stratified symplectic space given in
\cite{ls}.

The results announced in this paper are contained in \cite{bp2}.

In a subsequent paper, we will consider these spaces (as we have
already done in the simple case) in the complex setting and
therefore view them as a natural extension of the notion of toric
variety (cf. Remark~\ref{complex}).

For the definition of symplectic quasifold, quasitorus and every
related notion we refer the reader to \cite{p}.

\section{Stratifications by quasifolds}
We define the notion of space stratified by quasifolds in the
generality we need for our purposes. For the general definition of
stratification see \cite{GMcP1,GMcP2}.
\begin{defn}{\rm
Let $M$ be a compact topological space.
A {\it decomposition of $M$ by quasifolds} is a collection of
disjoint locally closed connected
quasifolds ${\cal T}_{F}$ ($F\in{\cal F}$),
called {\it pieces}, such that
\begin{enumerate}
\item The set $\cal F$ is finite and partially ordered.
\item $M=\bigcup_F{\cal T}_F$;
\item ${\cal T}_F\cap{\overline{\cal T}}_{F'}\neq\emptyset$ iff ${\cal
T}_F\subseteq{\overline {\cal T}}_{F'}$ iff $F\leq F'$.
\end{enumerate}}\end{defn}
We also require that ${\cal F}$ has a maximal element $F$ and that
the corresponding piece ${\cal T}_F$ is open and dense in $M$. We
call this piece the regular piece, the other pieces are called
singular. We will then say that $M$ is an {\em $n$-dimensional
compact space decomposed by quasifolds}, with $n$ the dimension of
the regular set.
\begin{remark}{\rm
A standard construction that is useful for the definition of
stratification is that of a cone over a compact space $L$
decomposed by quasifolds. We will call {\it cone over $L$},
denoted by $\stackrel{\circ}{C}(L)$, the space $[0,1)\times
L/\sim$, where two points $(t,l)$ and $(t',l')$ in $[0,1)\times L$
are equivalent if and only if $t=t'=0$. This space is itself a
space decomposed by quasifolds: for example when $L$ is a compact
quasifold the space $\stackrel{\circ}{C}(L)$ decomposes into two
pieces: one is the cone point, the other is given by the quasifold
$(0,1)\times L$. In fact we shall be considering a slightly more
complicated situation: let $t$ be a point in a quasifold $\cal T$,
B an open neighborhood of $t$ and $L$ a compact  space decomposed
by quasifolds. The decomposition of $L$ induces a decomposition of
the product $B\times\cono(L)$: to each piece $\cal L$ of $L$ there
corresponds the piece $B\times(0,1)\times{\cal L}$; to cover the
whole of $B\times\cono(L)$ we add a minimal piece, lying in the
closure of all other pieces, given by $B$ times the cone
point.}\end{remark} A stratification is a decomposition that is
locally well behaved.
\begin{defn}\label{stratificazione}{\rm Let $M$ be an $n$-dimensional
compact  space decomposed by quasifolds, the decomposition of  $M$
is said to be a {\it stratification by quasifolds} if each
singular piece ${\cal T}$, called {\it stratum}, satisfies the
following conditions:
\begin{enumerate}
\item[(i)] let $r$ be the dimension of $\cal T$,
for every point $t\in{\cal T}$ there exist an open neighborhood
$U$ of $t$ in $M$, an open neighborhood $B$ of $t$ in ${\cal T}$,
an $(n-r-1)$-dimensional compact  space $L$ decomposed by
quasifolds, called the {\it link} of $t$, and a homeomorphism
\linebreak $h\,\colon\,B\times \cono(L)\lorw U$ that preserves the
decompositions and that takes each piece of $B\times\cono(L)$
homeomorphically into the corresponding piece of $U$;
\item[(ii)] the decomposition of $L$ satisfies condition (i).
\end{enumerate}
}\end{defn}
The definition is recursive and, since the dimension of $L$ decreases
at each step, we end up, after a finite number of steps,
with links that are compact quasifolds.

\section{The construction}
Let $\d$ be a real vector space of dimension $n$, and let $\Delta$
be a convex polytope of dimension $n$ in the dual space $\ddu$. We
want to associate to the polytope $\D$ a family of compact spaces
that are suitably stratified by symplectic quasifolds. We
construct these spaces as symplectic quotients, following the
procedure which was first introduced by Delzant in \cite{delzant}.
Write the polytope as
\begin{equation}\label{polydecomp}
\D=\bigcap_{j=1}^d\{\;\mu\in\ddu\;|\;\langle\mu,X_j\rangle\geq\lambda_j\;\}
\end{equation}
for some elements $\xd$ in the vector space $\d$ and some real
numbers $\ld$. Let $Q$ be a quasilattice in the space $\d$
containing the elements $X_j$ (for example the one that is
generated by these elements) and let $\{\ed\}$ denote the standard
basis of $\rd$; consider the surjective linear mapping
$$
\begin{array}{cccc}\label{pi}
\pi \,\colon\,& \R^d & \lorw & \d\\
    &   e_j& \longmapsto & X_j.
\end{array}
$$
Consider the $n$-dimensional quasitorus $\d/Q$. The mapping $\pi$
induces a group homomorphism,
\begin{equation}\label{defdienne}
\Pi
\,\colon\, T^d=\rd/\Z^d\lorw \d/Q.
\end{equation}
We define $N$ to be the kernel of the mapping $\Pi$.
The mapping $\Pi$ defines an isomorphism
\begin{equation}\label{qtiso}
T^d/N\longrightarrow \d/Q.
\end{equation}
We construct a moment mapping for the Hamiltonian action of $N$ on
$\C^d$.
 Consider the mapping $J(\vz)=\sum_{j=1}^d
(|z_j|^2+\lambda_j)e_j^*$, where the $\lambda_j$'s are given in
(\ref{polydecomp}) and are uniquely determined by our choice of
normal vectors. The mapping $J$ is a moment mapping for the
standard action of $T^d$ on $\C^d$. Consider now the subgroup
$N\subset T^d$ and the corresponding inclusion of Lie algebras
$\iota\,\colon\,\n\rightarrow\R^d$. The mapping
$\Psi\,\colon\,\C^d\rightarrow \n^*$ given by $\Psi={\iota}^*\circ
J$ is a moment mapping for the induced action of $N$ on $\cd$. We
want to prove that the quotient $M=\zset/N$, endowed with the
quotient topology, is a space stratified by  quasifolds. Notice
that, by (\ref{defdienne}), the group $N$ is not necessarily
closed in $T^d$, moreover to each $\d$ there corresponds a whole
family of quotients, given by all possible choices of normal
vectors and of quasilattices $Q$ containing these vectors.

In our general setting, in which the polytope can be nonsimple,
the zero set $\zset$ is not in general a smooth submanifold of
$\R^{2d}$. Nonsimpleness of the polytope is responsible, like in
the rational case, for the decomposition in strata of the
quotient, whilst nonrationality produces the quasifold structure
of the strata.

To define the decomposition of $M$ in pieces we start by giving
some further definitions on the polytope $\D$.

Let us consider the open faces of $\Delta$. They can be described
as follows. For each such face $F$ there exists a possibly empty
subset $I_F\subset\{1,\ldots,d\}$ such that
\begin{equation}\label{facce}
F=\{\,\mu\in\Delta\;|\;\langle\mu,X_j\rangle=\lambda_j\;
\hbox{ if and only if}\; j\in I_F\,\}.
\end{equation}
A partial order on the set of all faces of $\D$ is defined by
setting $F\leq F'$ (we say $F$ contained in $F'$) if
$F\subseteq\overline{F'}$. The polytope $\Delta$ is the disjoint
union of its faces. Let $r_F=\hbox{card}(I_F)$; we have the
following definitions:

\begin{defn}{\rm A $p$-dimensional face $F$ of the polytope is
said to be {\em nonsimple} or {\em singular} if $r_F>n-p$.}
\end{defn}

\begin{defn}{\rm A $p$-dimensional face $F$ of the polytope is
said to be {\em simple} or {\em regular} if $r_F=n-p$.}
\end{defn}

\begin{prop}\label{phi} The $n$-dimensional quasitorus $D=\d/Q$ acts
continuously on the topological space $M=\zset/N$. Moreover $M$ is
compact and a continuous mapping $\Phi\,\colon\,M\lorw\d^*$ is
defined such that $\Phi(M)=\D$.\end{prop} \proof Consider the
exact sequence
\begin{equation}
\label{exactsequence}
0\lorw\d^*\stackrel{\pi^*}{\lorw}(\R^d)^*\stackrel{{\iota}^*}{\lorw}(\n)^*\lorw0.
\end{equation}
By (\ref{exactsequence}) we have that
the mapping $(\pi^*)^{-1}\circ J$ gives a well defined mapping on the quotient $M$,
we call this mapping $\Phi$. Moreover $\vz\in\zset$ if and only if
\begin{equation}
\label{formula} |z_j|^2=\langle \Phi(\vz),X_j\rangle-\lambda_j,\quad\quad j=1,\ldots, d.
\end{equation}
This implies that $\Phi(M)=\D$. Moreover properness of $J$ implies
that $M$ is compact.

\begin{remark}\label{stratidizset}{\rm From the proof of
Proposition~\ref{phi} we deduce that for any face $F$ of $\D$ the
set $\Phi^{-1}(F)$ is non-empty and is precisely given by
$\zset\cap\C^d_F/N$, where
$\cdf=\{\,\zd\in\cd\;|\;z_j=0\;\;\hbox{iff}\;\; j\in I_F\,\}$ (for
further detail cf. \cite{g,p}).}\end{remark}

\section{The stratification}
We are now ready to define the decomposition of $M$: the pieces
are given by ${\cal T}_F=\Phi^{-1}(F)$ with $F$ singular face, and
by the union, over all nonsingular faces $F$, of the sets
$\Phi^{-1}(F)$. This union is the  regular piece of the
decomposition while the ${\cal T}_F$'s are the singular pieces.
Let $\hbox{Int}(\Delta)$ be the open face of $\Delta$, then the
regular set contains $\Phi^{-1}(\hbox{Int}(\Delta))$. We will
label the regular set by the index $\hbox{Int}(\Delta)$ and call
it in short ${\cal T}_{\D}$. Remark~\ref{stratidizset} allows us
to characterize the pieces ${\cal T}_F$ in the standard way, by
the isotropy group attached to each of them.
\begin{remark}{\rm Let $F$ be a $p$-dimensional face. Let
${\mbox{\frak{s}}}^{F}=\{\,(y_1,\cdots,y_d)\;|\;y_j=0\;\;\hbox{if}\;\;
j\notin I_F\,\}$
and
$S^{F}=\{\,(Y_1,\cdots,Y_d)\in T^d\;|\;Y_j=1\;\;\hbox{if}\;\;j\notin I_F\,\}$.
The torus $S^F$ is the stabilizer of $T^d$ at any point
$\zd\in\zset\cap\C^d_F$ and ${\mbox{\frak{s}}}^{F}$ is its Lie
algebra. The stabilizer of $N$ on $\zset\cap\C^d_F$ is then the
$(r_F-n+p)$-dimensional subgroup $N^F$ of $N$ given by $N\cap
S^F$. Its Lie algebra, $\n^F$, is given by $\n\cap{\mbox{\frak{s}}}^{F}$.
Notice that the regular set $\zset_{\D}$, given by
the union, over all non-singular faces, of the sets
$\zset\cap\C^d_F$, has discrete stabilizer.}
\end{remark}
\begin{thm}\label{stratiquasifold}
The subset ${\cal T}_F$ of $M$ corresponding to each $p$-dimensional
singular face of $\D$ is a $2p$-dimensional quasifold.
In particular ${\cal T}_{\D}$
is a $2n$-dimensional quasifold. These subsets give a
decomposition by quasifolds of $M$.
\end{thm}
\begin{remark}{\rm
A singular face has at most dimension $n-2$, therefore a singular
piece has at most dimension $2n-4$}\end{remark}
\begin{remark}\label{strutturalisciaglobale}{\rm The decomposition
of $M$ is induced by the decomposition of $\zset$ given by the
manifolds $\zset\cap\C^d_F$, with $F$ singular, and the open
subset $\zset_{\D}$ of $\zset$. The quasifold structure of each
piece ${\cal T}_F$ is naturally induced by the smooth structure of
$\zset\cap\C^d_F$, while the quasifold structure of ${\cal
T}_{\D}$ is induced be the smooth structure of
$\zset_{\D}$}.\end{remark} Let $p\,\colon\,\zset\lorw M$ be the
projection, we have the following
\begin{thm}\label{stratisimplettici} Each piece ${\cal T}_F$ (${\cal T}_{\D}$)
of the decomposition of $M$ has a natural symplectic structure
induced by the quotient procedure, that is, its pull-back via $p$
coincides with the restriction of the standard symplectic form of
$\C^d$ to the manifold $\zset\cap\C^d_F$ ($\zset_{\D}$).\end{thm}
\begin{thm}\label{momentmapping} The restriction of the $D$-action
and of the mapping $\Phi$ to each piece of the space $M$ is
smooth, the action of $D$ is Hamiltonian and a moment mapping is
given by the restriction of $\Phi$.\end{thm}

Now we need to prove that our decomposition has a good local
behavior. Let $t$ be a point in the singular $2p$-dimensional
piece ${\cal T}_F$; we want to construct a link of $t$ satisfying
Definition~\ref{stratificazione}. Let
${\mbox{\frak{s}}}^{F}_{\SC}={\mbox{\frak{s}}}^{F}+i{\mbox{\frak{s}}}^{F}$
be the complexification of the Lie algebra
${\mbox{\frak{s}}}^{F}$. The mapping $J$ restricted to
${\mbox{\frak{s}}}^{F}_{\SC}$ gives rise to a moment mapping for
the action of $S^F$ on  ${\mbox{\frak{s}}}^{F}_{\SC}$, we denote
this mapping by $J_F\,\colon\,{\mbox{\frak{s}}}^{F}_{\SC}\lorw
\left({{\mbox{\frak{s}}}^{F}}\right)^*$. Consider now the
Hamiltonian action of the $(r_{F}-n+p)$-dimensional group $N^{F}$
on ${\mbox{\frak{s}}}^{F}_{\SC}$, induced by that of $S^F$: a
moment mapping is then given by $\psi_F=\iota_F^*\circ J_F$, where
$\iota_{F}\,\colon\,\n^{F}\lorw{{\mbox{\frak{s}}}^{F}}$ is the
inclusion map. In fact, using (\ref{facce}), it turns out that
$\psi_{F}=\sum_{j\in I_F}|z_j|^2\iota_{F}^*(e_j)$, hence
$\psi^{-1}_F(0)$ is a cone. We can now construct the link of $t$:
to begin with take the sphere ${\cal S}_{F,\epsilon}$ in
${\mbox{\frak{s}}}^{F}_{\SC}$ of radius $\epsilon$, centered in
$0$; we have that for any given $\epsilon>0$ the space
$\psi_F^{-1}(0)\cap {\cal S}_{F,\epsilon}$ is nonempty and is
acted on by the group $N^F$. Let us denote the quotient,
$\left(\psi_F^{-1}(0)\cap {\cal S}_{F,\epsilon}\right)/N^F$ by
$L_{F,\epsilon}$. We have the following:
\begin{lemma}\label{lemma} Let $t$ be a point in the singular piece
${\cal T}_F$. Then we can choose suitable open neighborhoods $B_t$
and $U_t$ of $t$, in ${\cal T}_F$ and $M$ respectively, and
an $\epsilon>0$, such that a decomposition preserving
homeomorphism
$h_{F}\,\colon\,B_t\times\cono(L_{F,\epsilon})\lorw U_t$ is defined.
The mapping $h_F$ restricted to each piece  is a
homeomorphism.
Moreover
$L_{F,\epsilon}$ is a
$(2n-2p-1)$-dimensional compact space decomposed by quasifolds.
\end{lemma}
\begin{thm}\label{locale}
Let $F$ be a singular face of the convex polytope $\Delta$ and $t$
be a point in ${\cal T}_F$. The compact space $L_{F,\epsilon}$ is
a link of $t$.
\end{thm}

The proof of Lemma~\ref{lemma} is based on
Theorem~\ref{stratiquasifold} and on the explicit construction of
the homeomorphism $h_{F}$ from $B_t\times\cono(L_{F,\epsilon})$
onto $U_{t}$. The proof of Theorem~\ref{locale} also consists in
exhibiting explicitly, at each step of the recursive definition of
link, a link for the point in consideration together with the
corresponding homeomorphism.
\begin{remark}\label{strutture}{\rm
Theorem~\ref{locale} proves that the decomposition of
$M$ is in fact a stratification,  a notion which is purely topological.
But, from
Theorems~\ref{stratiquasifold},$\,$\ref{stratisimplettici}, we know that
each piece of the stratification of $M$ has the structure of a
symplectic quasifold, naturally induced by that of $\C^d$.
We will call $M$ a {\em space stratified by symplectic quasifolds}.}
\end{remark}
\begin{remark}{\rm In the light of Theorem~\ref{momentmapping},
we can view the mapping $\Phi$ as a moment mapping for the action of
the $n$-dimensional
quasi-torus $D$ on the $2n$-dimensional compact space $M$ stratified
by symplectic quasifolds. By Prop~\ref{phi}, the image $\Phi(M)$ of the moment
mapping $\Phi$, is exactly
the polytope $\D$.}\end{remark}
\begin{remark}\label{complex}{\rm The remark above
emphasizes the relationship between the space $M$ and the polytope
$\D$, which is very neat in the symplectic setting. From the
complex point of view we have a compact space $X$, homeomorphic to
$M$, stratified by complex quasifolds; $X$ is $n$-dimensional and
is acted on by the complexified torus $\Dc$ of same dimension.
Such an action has a dense open orbit, corresponding to the open
set $\Phi^{-1}(\hbox{Int}(\D))$. }\end{remark}

\noindent \small{\sc Dipartimento di Matematica Applicata ``G. Sansone'',
Via S. Marta 3, 50139 Firenze, ITALY,  {\tt mailto:fiamma@dma.unifi.it}\\
                        and\\
                        Laboratoire Dieudonn\'e,
                        Universit\'e de Nice, Parc Valrose, 06108 Nice
                        Cedex 2, FRANCE, {\tt mailto:elisa@alum.mit.edu}}


\begin{thebibliography}{AAAAA}
\bibitem[BBFK1]{bbfk} G. Barthel, J.-P. Brasselet, K.-H. Fieseler,
L. Kaup, Equivariant intersection cohomology of toric varieties,
Algebraic geometry: Hirzebruch 70 (Warsaw, 1998), 45--68, {\em
Contemp. Math.} \textbf{241} Amer. Math. Soc., Providence, RI, 1999.

\bibitem[BBFK2]{bbfk2} G. Barthel, J.-P. Brasselet, K.-H. Fieseler,
L. Kaup, Combinatorial intersection cohomology for fans, {\em
arXiv:math.AG/0002181}.

\bibitem[BP1]{cx} F. Battaglia, E. Prato, Generalized toric varieties for
simple nonrational convex polytopes, {\em Intern. Math. Res.
Notices} \textbf{ 24} (2001), 1315-1337.

\bibitem[BP2]{bp2}  F. Battaglia, E. Prato, A symplectic realization of
convex polytopes that are neither rational nor simple, {\em in
preparation}.

\bibitem[BL]{bl} P. Bressler, V. Lunts, Intersection cohomology on
nonrational polytopes, {\em arXiv:math.AG/0002006}.

\bibitem[D]{delzant}  T. Delzant, Hamiltoniens p\'eriodiques et image convexe
de l'application moment, {\em Bull. Soc. Math. France} \textbf{ 116} (1988),
315--339.

\bibitem[GM1]{GMcP1} M. Goresky, R.MacPherson, Intersection homology II,
{\em Invent. Math.} \textbf{ 71} (1983) 77--129.

\bibitem[GM2]{GMcP2} M. Goresky, R.MacPherson, Stratified Morse Theory, Springer Verlag, New York, 1988.

\bibitem[G]{g} V. Guillemin, Moment maps and combinatorial invariants
of Hamiltonian $T^n$-spaces, Progress in Mathematics 122,
Birkh\"auser, Boston, 1994.

\bibitem[LS]{ls}  R. Sjamaar, E. Lerman, Stratified symplectic spaces and reduction,
{\em Ann. of Math.} \textbf{ 134} (1991), 375-422.

\bibitem[P]{p}  E. Prato, Simple non-rational convex polytopes via
symplectic geometry, {\em Topology} \textbf{ 40} (2001), 961-975.

\end{thebibliography}
\end{document}